\documentclass[reqno,11pt]{amsart}

\usepackage{amssymb}
\usepackage{mathtools}
\usepackage[all,cmtip]{xy}
\usepackage{leftidx}
\newtheorem{Proposition}{Proposition}[section]

\newtheorem{Definition}[Proposition]{Definition}
\newtheorem{Lemma}[Proposition]{Lemma}
\newtheorem{Conjecture}[Proposition]{Conjecture}
\newtheorem{Theorem}[Proposition]{Theorem}
\newtheorem{MainTheorem}{Theorem}

\DeclareMathOperator{\Val}{Val}

\newcommand{\C}{\mathbb{C}}
\newcommand{\CP}{\mathbb{CP}}
\newcommand{\R}{\mathbb{R}}

\begin{document}
\title{Unitarily invariant valuations and Tutte's sequence}
\author{Andreas Bernig}

\email{bernig@math.uni-frankfurt.de}

\address{Institut f\"ur Mathematik, Goethe-Universit\"at Frankfurt,
Robert-Mayer-Str. 10, 60629 Frankfurt, Germany}


\begin{abstract}
We prove Fu's power series conjecture which relates the algebra of isometry invariant valuations on complex space forms to a formal power series from combinatorics which was introduced by Tutte. The $n$-th coefficient of this series is the number of triangulations of a triangle with $3n$ internal edges; or the number of intervals in Tamari's lattice $Y_n$.   
\end{abstract}

\thanks{{\it MSC classification}:  53C65, 
05A15 
\\ Supported
 by DFG grant BE 2484/5-2}
\maketitle

\section{Statement of the result and background}

Thanks to the groundbreaking work of Alesker \cite{alesker_mcullenconj01, alesker04_product, alesker_survey07}, the space of valuations on manifolds (i.e. finitely additive functionals on some class of sufficiently regular sets) was endowed with a product structure which satisfies a version of Poincar\'e duality. 

The product of invariant valuations on isotropic manifolds (i.e. Riemannian manifolds such that the isometry group acts transitively on the sphere bundle) encodes the kinematic formulas on such spaces \cite{bernig_fu06, bernig_fu_solanes}. Using this powerful new branch of integral geometry, which is called \emph{algebraic integral geometry}, it was possible to write down in explicit form the kinematic formulas on all complex space forms \cite{bernig_fu_hig,bernig_fu_solanes}. We refer to \cite{bernig_aig10, fu_barcelona} for a survey on these developments.

The starting point is a theorem by J.~Fu which describes the algebra of invariant valuations on hermitian space $(\C^n,\mathbb U(n))$.
\begin{Theorem} \label{thm_fu}
The algebra $\Val^{\mathrm U(n)}$ of translation-invariant, continuous and $\mathrm U(n)$-invariant valuations on $\C^n$ is isomorphic to the polynomial algebra
\begin{displaymath}
 \R[t,s]/(f_{n+1},f_{n+2}),
\end{displaymath}
where $\deg t=1, \deg s=2$, and
\begin{displaymath}
\log(1+t+s)=\sum_{i \geq 1} f_i(t,s)
\end{displaymath}
is the decomposition into homogeneous components. 
\end{Theorem}

The valuation $t \in \Val^{\mathrm U(n)}$ equals, up to normalization, the first intrinsic volume. The valuation $s$ equals on a compact convex subset $K \subset \C^n$ the measure of complex hyperplanes intersecting $K \subset \C^n$ (again up to normalization).   

This theorem was the main entry into a deep study of the integral geometry of $(\C^n,\mathrm U(n))$. In \cite{bernig_fu_hig}, several geometrically interesting bases of $\Val^{\mathrm U(n)}$ were introduced, their mutual relations were described and the kinematic formulas were explicitly written down.   

After these results in the flat case, the next challenge was to find a similar approach in the curved case. Given a real number $\lambda$, we denote by $\mathbb{CP}^n_\lambda$ the complex space form of constant holomorphic curvature $4\lambda$, endowed with the group $G_\lambda$ of holomorphic isometries. If $\lambda>0$, then $\CP^n_\lambda$ is a rescaling of complex projective space. If $\lambda<0$, then $\CP^n_\lambda$ is a rescaling of complex hyperbolic space, while $\CP^n_0$ is the flat space $\C^n$. 

A natural guess is that the algebra structure of $\mathcal{V}(\CP^n_\lambda)^{G_\lambda}$ can be described in a way analogous to Theorem \ref{thm_fu}, with some modification of the $f_i$ depending on $\lambda$. Based on numerical evidence, J.~Fu stated in several talks around 2008 and in written form in  \cite[Conjecture 6.14]{fu_barcelona} the following conjecture. 

\begin{Conjecture}[Fu's power series conjecture] \label{conj_fu_conjecture}
Define Tutte's series as
\begin{displaymath}
\tau(\lambda):=\sum_{i=1}^\infty \frac{2(4i+1)!}{(i+1)!(3i+2)!}\lambda^i=\lambda+3\lambda^2+13\lambda^3+68\lambda^4+399\lambda^5+\ldots
\end{displaymath}
Then the algebra of invariant valuations on $\mathbb{CP}_\lambda^n$ is isomorphic to 
\begin{displaymath}
 \R[[t,s]]/(\bar f_{n+1}^\lambda,\bar f_{n+2}^\lambda), 
\end{displaymath}
where the formal power series $\bar f_k^\lambda(t,s) \in \R[[t,s]]$ is defined as the degree $k$-part in the expansion of
\begin{displaymath}
\log\left(1+t+s+\tau(\lambda)\right).
\end{displaymath}
Here $t$ is of degree $1$, $s$ is of degree $2$, and $\lambda$ is of degree $(-2)$. 
\end{Conjecture}

In this conjecture, $t$ denotes a certain multiple of the first intrinsic volume (which may be defined canonically on any Riemannian manifold), and $s$ is the average Euler characteristic of the intersection with a totally geodesic complex hyperplane in $\CP^n_\lambda$. 

If we prefer an ideal in $\R[t,s]$ instead of $\R[[t,s]]$, we may truncate $\bar f_{n+1}^\lambda,\bar f_{n+2}^\lambda$ and add the monomials $t^{2n+1},\ldots,s^nt$. This is the original form of the conjecture, see \cite[Conjecture 1]{bernig_fu_solanes}. It will also come out of the proof that all $\bar f_i^\lambda, i>n$, vanish on $\CP_\lambda^n$, as conjectured in \cite[Conjecture 6.14]{fu_barcelona}  

The first few terms of the conjecture can be confirmed using the template method and some formulas involving binomial coefficients. These formulas can be shown using Zeilberger's algorithm, but get more and more involved for higher orders of $\lambda$, compare the discussion in \cite{fu_barcelona}. 

Interestingly, the sequence of coefficients $1,3,13,68,399,2530,\ldots$ has various combinatorial interpretations. Tutte  \cite{tutte62} has shown that the coefficient of $\lambda^i,i>0$ in $\tau(\lambda)$ is the number of non-isomorphic planar triangulations of a triangle with $3i$ internal edges. In \cite[Section 4]{cori_schaeffer}, the same coefficient appears as the number of description trees of type $(1,1)$ and size $i$. In \cite{chapoton} it was shown that this coefficient also equals the number of intervals in Tamari's lattice $Y_i$. In \cite{bernardi_bonichon} some explicit bijections explainig these numerical coincidences are constructed. See \cite{sloane} for many other appearances of this sequence.

Some years after the statement of the conjecture, and without using it, the integral geometry of complex space forms was worked out in \cite{bernig_fu_solanes}. A surprising result, which is based on some computations and not yet fully understood, is the following. 
 
\begin{Theorem} \label{thm_isomorphism}
The algebras $\Val^{\mathrm U(n)}$ and $\mathcal{V}(\CP^n_\lambda)^{G_\lambda}$ are isomorphic for each $\lambda$. 
\end{Theorem} 

There are many different isomorphisms, and one of them, given in Proposition \ref{prop_isomorphism} below, will play a central role in this paper.

Despite the progress in integral geometry of complex space forms (see also \cite{abardia_gallego_solanes, alesker03_un,  bernig_fu_solanes_proceedings, wannerer_area_measures, wannerer_unitary_module}), the original conjecture remained previously open. In this paper we are going to prove it. 

\begin{MainTheorem} \label{mainthm_fu_conjecture}
Fu's power series conjecture is true.
\end{MainTheorem}

Our proof is a mixture of the template method and some algebraic manipulations of generating functions. One of the main ingredients is the fact that $\tau(\lambda)$ satisfies an algebraic equation over $\R(\lambda)$, which implies that some auxiliary power series in $s$ and $\lambda$, which comes out of the template method, is of non-positive degree (where $\deg s=2, \deg \lambda=-2$). This part of the proof uses the theory of holonomic functions. We refer to the very recent lecture notes \cite{sattelberger_sturmfels} for more information on this subject.  

We do not know if there is a more direct link between integral geometry and the combinatorics of triangulations, Tamari's lattice or description trees. It is a priori not even clear that the coefficients in the power series in Fu's conjecture have to be  integers and not just any real numbers. Also the way they appear in the conjecture (via the relations satisfied in some algebra) does not seem to be related to any counting of objects.   
 
\subsection*{Thanks}
I would like to thank Joseph Fu for many interesting discussions about hermitian integral geometry and for the (still ongoing) collaboration on this subject. I thank Anna-Laura Sattelberger for her helpful explanations about D-modules and holonomic functions. I thank the anonymous referee for useful suggestions for improving the presentation.

\section{An isomorphism}

\begin{Proposition} \label{prop_isomorphism}
The map given by $t \mapsto t, s \mapsto \frac{s}{1-\lambda s}$ induces an algebra isomorphism between $\Val^{\mathrm U(n)}$ and $\mathcal{V}(\CP^n_\lambda)^{G_\lambda}$. 
\end{Proposition}

\proof
For each $\lambda \in \R$, the map $\R[[t,s]] \to \R[[t,s]], t \mapsto \frac{t}{\sqrt{1-\lambda s}}, s \mapsto \frac{s}{1-\lambda s}$ covers an isomorphism $\Val^{\mathrm U(n)} \to \Val^{\mathrm U(n)}$, since the ideal defining $\Val^{\mathrm U(n)}$ is homogeneous (see Theorem \ref{thm_fu}). By \cite[Theorem 3.17]{bernig_fu_solanes}, the map $t \mapsto t \sqrt{1-\lambda s}, s \mapsto s$ covers an isomorphism from $\Val^{\mathrm U(n)}$ to $\mathcal V(\CP_\lambda^n)^{G_\lambda}$. The displayed map is the composition of these two isomorphisms, and hence an isomorphism as well.
\endproof

The proposition implies that the map $t \mapsto t, s \mapsto \frac{s}{1+\lambda s}$ is an isomorphism from $\mathcal{V}(\mathbb{CP}_\lambda^n)^{G_\lambda}$ to $\Val^{\mathrm U(n)}$. Therefore, Theorem \ref{mainthm_fu_conjecture} is equivalent to saying that $\Val^{\mathrm U(n)}$ is isomorphic to $\mathbb \R[[t,s]]/(f_{n+1}^\lambda,f_{n+2}^\lambda)$, where $f_k^\lambda(t,s)$ is the degree $k$-component in 
\begin{displaymath}
\log\left(1+t+\frac{s}{1+\lambda s}+\tau(\lambda)\right).
\end{displaymath}
We will prove it in this form.

\section{The template method}

We first need some preparations. All computations will be done with formal power series, so no convergence is required.  
 
Let us recall the formal power series 
\begin{align*}
 \log(1+x) & = \sum_{k=1}^\infty \frac{(-1)^k}{k}x^k,\\
 (1+x)^\alpha & = \sum_{k=0}^\infty \binom{\alpha}{k} x^k.
\end{align*}

Taking $\alpha:=-\frac12$ and rewriting the binomial coefficient we get
\begin{equation} \label{eq_binomial_series}
  \frac{1}{\sqrt{1-4x}}=\sum_{l=0}^\infty \binom{2l}{l} x^l.
 \end{equation}

\begin{Lemma}
For all $m \geq 0$ we have 
\begin{align} \label{eq_combinatorial_series} 
\sum_{\substack{k \equiv m \mod 2\\k>0}} \binom{k+m}{\frac{k+m}{2}} \frac{x^k}{k} & =Q_1^m\left(\frac{1}{x}\right)+Q_2^m \left(\frac{1}{x}\right)\sqrt{1-4x^2} \nonumber\\
& \quad \underbrace{-\binom{m}{\frac{m}{2}} \log(1+\sqrt{1-4x^2})}_{\text{ if } m \text{ is even}},
\end{align}
where $Q_1^m,Q_2^m$ are polynomials of degree $m$. If $m$ is even, then $Q_1^m,Q_2^m$ are even polynomials, and $Q_2^m$ does not contain an absolute term. If $m$ is odd, then $Q_1^m,Q_2^m$ are odd polynomials, and the logarithmic term does not appear. 
\end{Lemma}

\proof
Set
\begin{displaymath}
 Q_1^m(y) :=  \sum_{\substack{1 \leq i \leq m\\i \equiv m \mod 2}}  \binom{m-i}{\frac{m-i}{2}} \frac{y^i}{i}
\end{displaymath}
and 
\begin{displaymath}
F(x):=\sum_{\substack{k \equiv m \mod 2\\ k>0}} \binom{k+m}{\frac{k+m}{2}} \frac{x^k}{k}-Q_1^m\left(\frac{1}{x}\right)+ \underbrace{\binom{m}{\frac{m}{2}} \log(1+\sqrt{1-4x^2})}_{\text{ if } m \text{ is even}}.
\end{displaymath}

Using \eqref{eq_binomial_series} we compute
\begin{align*}
F'(x) & =\sum_{\substack{k \equiv m \mod 2\\ k>0}} \binom{k+m}{\frac{k+m}{2}} x^{k-1}+\sum_{\substack{1 \leq i \leq m\\i \equiv m \mod 2}}  \binom{m-i}{\frac{m-i}{2}} \frac{1}{x^{i+1}}-\underbrace{\binom{m}{\frac{m}{2}} \frac{1-\sqrt{1-4x^2}}{x \sqrt{1-4x^2}}}_{\text{ if } m \text{ is even}}\\
& = \sum_{j>\frac{m}{2}} \binom{2j}{j} \frac{x^{2j}}{x^{m+1}}+\sum_{0 \leq j<\frac{m}{2}} \binom{2j}{j} \frac{x^{2j}}{x^{m+1}}+\underbrace{\binom{m}{\frac{m}{2}} \frac{x^m}{x^{m+1}}-\binom{m}{\frac{m}{2}} \frac{1}{x \sqrt{1-4x^2}}}_{\text{ if } m \text{ is even}}\\
& = \frac{1}{x^{m+1}\sqrt{1-4x^2}}- \underbrace{\binom{m}{\frac{m}{2}} \frac{1}{x\sqrt{1-4x^2}}}_{\text{ if } m \text{ is even}}.
\end{align*}

We now look for some polynomial $Q_2^m$ such that 
\begin{displaymath}
 \left[ Q_2^m \left(\frac{1}{x}\right)\sqrt{1-4x^2}  \right]'=F'(x)=\frac{1}{x^{m+1}\sqrt{1-4x^2}}- \underbrace{\binom{m}{\frac{m}{2}}  \frac{1}{x\sqrt{1-4x^2}}}_{\text{ if } m \text{ is even}}.
\end{displaymath}

This is equivalent to 
\begin{displaymath}
(4-y^2) \frac{d}{dy} Q_2^{m}(y)-\frac{4}{y} Q_2^m(y)=y^{m+1}-\underbrace{\binom{m}{\frac{m}{2}} y}_{\text{ if } m \text{ is even}}
\end{displaymath}
and admits the solution 
\begin{displaymath}
 Q_2^m(y):=- \sum_{\substack{1 \leq i \leq m\\i \equiv m \mod 2}}  \frac{2^{m-i}i!!(m-1)!!}{i (i-1)!!m!!} y^i,
\end{displaymath}
where we use the standard notation $i!!=i \cdot (i-2) \cdot (i-4) \cdot \ldots$, with the convention $(-1)!!=0, 0!!=1$.

It follows that with these choices of $Q_1^m,Q_2^m$, the derivative of \eqref{eq_combinatorial_series} with respect to $x$ vanishes, hence the equation is correct up to a constant. If $m$ is odd, then both sides of the equation are odd, hence the constant must be $0$. If $m$ is even, we may add the constant to $Q_1^m$.  
\endproof

\begin{Definition}
A formal power series $b=b(s,\lambda)$ is called \emph{of degree $\leq 2m$} if each monomial $s^i \lambda^j$ appearing in $b$ has $2i-2j \leq 2m$.  
Equivalently, 
\begin{displaymath}
b^{(k)}(0):= \left(\left.\frac{\partial}{\partial \lambda}\right|_{\lambda=0}\right)^k b(s,\lambda)
\end{displaymath}
is a polynomial in $s$ of degree $\leq 2(m+k)$ for each $k \geq 0$. 
\end{Definition}

It is clear that the sum of two series of degrees $\leq 2m$ is again of degree $\leq 2m$ and that the product of a series of degree $\leq 2m_1$ and a series of degree $\leq 2m_2$ is a series of degree $\leq 2(m_1+m_2)$. 

If $c>0$ and $b$ has no constant term, then 
\begin{displaymath}
\log(c+b)=\log(c)+\log\left(1+\frac{b}{c}\right)=\log(c)+\sum_{n=1}^\infty \frac{(-1)^{n+1}}{n} \frac{b^n}{c^n}
\end{displaymath}
is a well-defined formal power series. If moreover $b$ is of degree $\leq 0$, then $\log(c+b)$ is also of degree $\leq 0$. 

\begin{Proposition} \label{prop_chapoton_dominant}
Let $r:=\frac{s}{1+\lambda s}+\tau(\lambda)$, where $\tau(\lambda)$ is Tutte's sequence. Then the formal power series
 \begin{displaymath}
b(s,\lambda):=  s+(1+\lambda s) \sqrt{(1+r)^2-4s}
\end{displaymath}
 is of degree $\leq 0$.
\end{Proposition}

We postpone the technical proof to the next section. 

Given a formal power series $p$ in $t,s$, let us write 
\begin{displaymath}
\int_{\C^n} p:=p_{2n}(D_n^\C),
\end{displaymath} 
where $p_{2n}$ denotes the $(2n)$-homogeneous component of $p$, considered as an element of $\Val_{2n}^{\mathrm U(n)}$ under the isomorphism from Theorem \ref{thm_fu}, and $D_n^\C \subset \C^n$ is the unit ball.  By \cite[Equation (134)]{fu_barcelona} we have 
\begin{equation} \label{eq_fu_integral}
\int_{\C^n} s^k t^{2n-k}=\binom{2n-2k}{n-k}.
\end{equation}
This implies that for every $p$ and $n \geq 0$ we have
\begin{equation} \label{eq_mult_s}
 \int_{\C^{n}} sp=\int_{\C^{n-1}} p,
\end{equation}
where we formally set $\int_{\C^{n-1}}p:=0$ if $n=0$.

Let us write $p \equiv q$ if $p,q$ are formal power series in $t,s$ such that $\int_{\C^n} p=\int_{\C^n} q$ for all $n \geq 0$. From \eqref{eq_fu_integral} it follows that 
\begin{displaymath}
t^k \equiv \begin{cases} \binom{k}{\frac{k}{2}} s^{\frac{k}{2}} & \text{if } k \text{ is even,}\\
0 & \text{if } k \text{ is odd}.
\end{cases}
\end{displaymath}

If $p \equiv q$ then \eqref{eq_mult_s} implies that $sp \equiv sq$.

\begin{Lemma} \label{lemma_construction_h}
Let $r=\frac{s}{1+\lambda s}+\tau(\lambda)$ as above. There exists a formal power series $h_m(s,\lambda) \in \R[[s,\lambda]]$ of degree $\leq 2m$ such that  
\begin{equation} \label{eq_reduction_to_s}
 t^{m}\log\left(1+t+r\right) \equiv h_m(s,\lambda).
 \end{equation}
 \end{Lemma}

\proof
 We compute 
\begin{align*}
 t^{m} & \log(1+t+r) = \sum_{n=1}^\infty t^{m} \frac{(-1)^{n+1}}{n}(t+r)^n\\
 & = \sum_{n=1}^\infty \frac{(-1)^{n+1}}{n} \sum_{k=0}^n \binom{n}{k} t^{k+m} r^{n-k}\\
 & \equiv \sum_{n=1}^\infty \frac{(-1)^{n+1}}{n} \sum_{k \equiv m \mod 2} \binom{n}{k} \binom{k+m}{\frac{k+m}{2}} s^{\frac{k+m}{2}} r^{n-k}\\
 & =\underbrace{ \binom{m}{\frac{m}{2}}s^{\frac{m}{2}} \log(1+r)}_{\text{ if } m \text{ even}}+\sum_{k \equiv m \mod 2, k>0} \binom{k+m}{\frac{k+m}{2}} s^{\frac{k+m}{2}} \sum_{n=k}^\infty \frac{(-1)^{n+1}}{n}  \binom{n}{k}  r^{n-k}\\
 & = \underbrace{ \binom{m}{\frac{m}{2}}s^{\frac{m}{2}} \log(1+r)}_{\text{ if } m \text{ even}}+(-1)^{m+1} \sum_{k \equiv m \mod 2, k>0} \binom{k+m}{\frac{k+m}{2}} \frac{s^{\frac{k+m}{2}}}{k(1+r)^{k}}\\
 & = h_m(s,\lambda), 
\end{align*}
where 
\begin{align*}
h_{m}(s,\lambda) & := \binom{m}{\frac{m}{2}}s^{\frac{m}{2}} \left(\log(1+r)+\log\left(1+\sqrt{1-\frac{4s}{(1+r)^2}}\right)\right)\\
 & \quad  +(-1)^{m+1} s^{\frac{m}{2}} \left(Q_1^{m}\left(\frac{1+r}{\sqrt{s}}\right)+Q_2^m\left(\frac{1+r}{\sqrt{s}}\right)\sqrt{1-\frac{4s}{(1+r)^2}}\right),                      
\end{align*}
and the logarithmic terms only appear if $m$ is even. 

It remains to show that $h_m$ is of degree $\leq 2m$. Consider a monomial in the $Q_1$-term. Let $0 \leq i \leq m, i \equiv m \mod 2$. Then 
\begin{align*}
s^{\frac{m}{2}} \left(\frac{1+r}{\sqrt{s}}\right)^{i} & = s^{\frac{m-i}{2}} \left(1+\frac{s}{1+\lambda s}+\tau(\lambda)\right)^{i}\\
& = s^{\frac{m-i}{2}} \sum_{a=0}^{i} \binom{i}{a}\left(\frac{s}{1+\lambda s}\right)^a (1+\tau(\lambda))^{i-a} \\
& = \sum_{a=0}^{i} \binom{i}{a} s^{\frac{m-i}{2}+a} \frac{(1+\tau(\lambda))^{i-a}}{(1+\lambda s)^a}.
\end{align*}
Since $a \leq i$, each summand is of degree $\leq 2m$.

Next, we consider a monomial in the $Q_2$-term. Let $1 \leq i \leq m, i \equiv m \mod 2$. Then 
\begin{align*}
s^\frac{m}{2} \left(\frac{1+r}{\sqrt{s}}\right)^{i} \sqrt{1-\frac{4s}{(1+r)^2}} & = s^\frac{m-1}{2} \left(\frac{1+r}{\sqrt{s}}\right)^{i-1} \sqrt{(1+r)^2-4s}\\
& =  s^\frac{m-1}{2} \left(\frac{1+r}{\sqrt{s}}\right)^{i-1} \frac{1}{1+\lambda s} (b(s,\lambda)-s).
\end{align*}

By Proposition \ref{prop_chapoton_dominant}, $(b(s,\lambda)-s)$ is of degree $\leq 2$. As we have seen above, $s^\frac{m-1}{2} \left(\frac{1+r}{\sqrt{s}}\right)^{i-1}$ is of degree $\leq 2(m-1)$. Hence the whole term is of degree $\leq 2m$.

Let us finally consider the logarithmic term, which only appears if $m$ is even. We have 
\begin{align*}
\log(1+r)+& \log\left(1+\sqrt{1-\frac{4s}{(1+r)^2}}\right) \\
& = \log\left(1+r+\sqrt{(1+r)^2-4s}\right)\\
  & = \log\left(1+\lambda s + (1+\lambda s)\tau(\lambda)+b(s,\lambda)\right)-\log(1+\lambda s).
\end{align*}

By Proposition \ref{prop_chapoton_dominant}, this term is of degree $\leq 0$.  Since we multiply it by $\binom{m}{\frac{m}{2}} s^{\frac{m}{2}}$, we get a term of degree $\leq m \leq 2m$. 
\endproof

\proof[Proof of Theorem \ref{mainthm_fu_conjecture}]
Let us write 
\begin{displaymath}
f_n^\lambda=\sum_{i=0}^\infty f_{n,i}(t,s)\lambda^i, \quad \deg f_{n,i}(t,s)=n+2i,
\end{displaymath}
and 
\begin{displaymath}
h_m(s,\lambda)=\sum_{i \geq 0} h_{m,i}(s) \lambda^i.
\end{displaymath}

Let us show that $\int_{\C^n} t^{m}s^l f_{n+1}^\lambda=0$ for all $n,m,l$ or equivalently that  
$\int_{\C^n} t^{m}s^l f_{n+1,i}=0$ for all $n,i,m,l$. Alesker duality \cite[Theorem 0.8]{alesker04_product} then implies that $f_{n+1,i}=0$ for all $i$ and hence $f_{n+1}^\lambda=0$ on $\C^n$.

By degree reasons, we have to show this only if $n=m+2l+2i+1$. Now
\begin{align*}
\int_{\C^{m+2l+2i+1}}  t^{m}s^l f_{m+2i+2,i}(t,s) & =\int_{\C^{m+2l+2i+1}}  t^{m}s^l \sum_q f_{q,i}(t,s)\\
& =\int_{\C^{m+2l+2i+1}}  s^lh_{m,i}(s),
\end{align*}
where the last equation follows from comparing the coefficients of $\lambda^i$ in  \eqref{eq_reduction_to_s}. 

By Lemma \ref{lemma_construction_h}, we have $\deg s^lh_{m,i} \leq 2l+2(m+i)<2(m+2l+2i+1)$, hence the integral on the right hand side vanishes. This finishes the proof that $f_{n+1}^\lambda=0$ on $\C^n$ for all $n$. If $i>1$, then $f_{n+i}^\lambda=f_{(n+i-1)+1}^\lambda=0$ on $\C^{n+i-1}$ by what we have shown and hence $f_{n+i}^\lambda=0$ on $\C^n$ by restriction from $\C^{n+i-1}$ to $\C^n$. 

Let us write 
\begin{displaymath}
\Phi:\R[[t,s]] \to \mathcal{V}(\CP_\lambda^n)^{G_\lambda}
\end{displaymath}
for the natural algebra morphism. This map is surjective, and we have shown that $\bar f_{n+1}^\lambda,\bar f_{n+2}^\lambda \in \ker \Phi$. To finish the proof, we have to show that these two power series generate the kernel. 

By \cite[Theorem 3.1.1.]{alesker_val_man4} the algebra $\mathcal{V}(\CP_\lambda^n)$ has a natural filtration 
\begin{displaymath}
\mathcal{V}(\CP_\lambda^n) = W_0(\CP_\lambda^n) \supset  W_1(\CP_\lambda^n) \supset \cdots \supset  W_{2n}(\CP_\lambda^n), 
\end{displaymath}
which is compatible with the product structure (where we put $W_k(\CP_\lambda^n)=\{0\}$ for $k>2n$). 

Fix a point $p_0 \in \CP^n_\lambda$. By \cite[Theorem 3.1.2]{alesker_val_man4} we have a linear map 
\begin{displaymath}
\Xi_k:W_k(\CP_\lambda^n)^{G_\lambda} \to \Val_k(T_{p_0}\CP^n_\lambda)^{\mathrm U(n)} \cong \Val_k^{\mathrm U(n)}.
\end{displaymath} 

Since $t \in  W_1(\CP_\lambda^n)^{G_\lambda}, s \in  W_2(\CP_\lambda^n)^{G_\lambda}$ it follows that if the lowest degree of $h \in \R[[t,s]]$ is $\geq k$, then $\Phi(h) \in W_k(\CP_\lambda^n)^{G_\lambda}$, and $\Xi_k \circ \Phi(h)$ equals the $k$-homogeneous part of $h$, which we denote by $h_k$. 

Suppose now that $h \in \R[[t,s]]$ is in the kernel of $\Phi$ and is of lowest degree $k$. Then $h_k=\Xi_k \circ \Phi(f)=0$ in $\Val^{\mathrm U(n)}$. Theorem \ref{thm_fu} implies the existence of a polynomial $p_k \in \R[t,s]$ of degree $(k-n-1)$ and a polynomial $q_k \in \R[t,s]$ of degree $(k-n-2)$ such that $h_k=p_k f_{n+1}+q_k f_{n+2}$. The lowest degree of $h-p_k \bar f_{n+1}^\lambda-q_k \bar f_{n+2}^\lambda \in \R[[t,s]]$ is then at least $(k+1)$. Repeating the same argument, we find sequences of polynomials $p_i,q_i \in \R[t,s],i=k,k+1,\ldots$ with $\deg p_i=i-n-1, \deg q_i=i-n-2$ such that $h-\sum_{i=k}^l (p_i \bar f_{n+1}^\lambda-q_i \bar f_{n+2}^\lambda)$ is of lowest degree at least $(l+1)$ for each $l \geq k$. Setting $p:=\sum_{i=k}^\infty p_i \in \R[[t,s]], q:=\sum_{i=k}^\infty q_i \in \R[[t,s]]$, we find that $h=p \bar f_{n+1}^\lambda+q \bar f_{n+2}^\lambda$.
\endproof

\section{Proof of Proposition \ref{prop_chapoton_dominant}}

The proof in this section is based on some terminology and results from D-modules and holonomic functions. We refer to \cite{coutinho, sattelberger_sturmfels, zeilberger} for more information on this subject. For the reader's convenience, we have spelled out the details so that it should be possible to understand the proof without any prior knowledge of holonomic functions. 

Let us first sketch the idea. Tutte has shown that $\tau(\lambda)$ is an algebraic function \cite[Equations (4.8), (4.9)]{tutte62} (see also \cite[Theorem 4]{cori_schaeffer}, \cite[Equation (10)]{chapoton}). The formal power series $b(s,\lambda)$ is therefore the composition of a holonomic and an algebraic function and is itself holonomic \cite[Proposition 2.3]{sattelberger_sturmfels}. The holonomic rank turns out to be $4$. In particular, if we consider $b$ as a formal power series in $\lambda$ with coefficients in the space of formal power series in $s$, it satisfies some differential equation of degree $4$ whose coefficients are polynomials in $s$ and $\lambda$. The differential equation gives a recursive relation for the higher derivatives $\left.\frac{d^i}{d\lambda^i}\right|_{\lambda=0} b(s,\lambda)$ which will be enough to prove the proposition.  

Let us now work out the details. 

Define an algebraic power series $\phi(\lambda)$ by the equation $P(\lambda,\phi(\lambda))=0$, where $P(\lambda,\phi):=\phi-\lambda(1+\phi)^4 \in \R[\lambda,\phi]$ is irreducible. Then, by \cite[Equation (10)]{chapoton}, $\tau(\lambda)=\phi(\lambda)(1-\phi(\lambda)-\phi(\lambda)^2)$.

\begin{Lemma}
 There exists some $P_0 \in \R(\lambda)[\phi]$ of degree $3$ such that 
 \begin{displaymath}
  \phi'(\lambda)=P_0(\lambda,\phi(\lambda)).
 \end{displaymath}
\end{Lemma}

\proof
Taking derivatives of the equation $P(\lambda,\phi(\lambda))=0$ yields 
\begin{displaymath}
\phi'(\lambda)=- \frac{\frac{\partial P}{\partial \lambda}(\lambda,\phi(\lambda))}{\frac{\partial P}{\partial \phi}(\lambda,\phi(\lambda))}=\frac{(1+\phi(\lambda))^4}{1-4\lambda(1+\phi(\lambda))^3}.
\end{displaymath}

Applying the extended Euclidean algorithm to $P$ and $\frac{\partial P}{\partial \phi}$ gives polynomials $U,V \in \R(\lambda)[\phi]$ such that $U P +V \frac{\partial P}{\partial \phi}=1$, and we deduce that 
\begin{displaymath}
\phi'(\lambda)=- \frac{V \frac{\partial P}{\partial \lambda}(\lambda,\phi(\lambda))}{V \frac{\partial P}{\partial \phi} (\lambda,\phi(\lambda))}=-V \frac{\partial P}{\partial \phi}(\lambda,\phi(\lambda)).
\end{displaymath}

Let $P_0 \in \R(\lambda)[\phi]$ be the remainder of $-V \frac{\partial P}{\partial \phi}$ by division by $P$. Then $\phi'(\lambda)=P_0(\lambda,\phi(\lambda))$. Doing the computations explicitly, we find that 
\begin{displaymath}
P_0(\lambda,\phi)=\frac{12 \lambda \phi^3+52 \lambda \phi^2+4\lambda \phi-36\lambda+9\phi}{(256\lambda-27)\lambda}.
\end{displaymath}
\endproof

\begin{Lemma}
Set 
\begin{displaymath}
F(s,\lambda,\phi):=\sqrt{\left[1+\lambda s+s+(1+\lambda s)\phi(1-\phi-\phi^2)\right]^2-4s(1+\lambda s)^2},
\end{displaymath}
so that $b(s,\lambda)=s+F(s,\lambda,\phi(\lambda))$. 
Then there exists $Q_1 \in \R(s,\lambda)[\phi]$ of degree $\leq 3$ with 
\begin{displaymath}
\frac{d}{d\lambda} F(s,\lambda,\phi(\lambda)) = Q_1(s,\lambda,\phi(\lambda)) F(s,\lambda,\phi(\lambda)).
\end{displaymath}
\end{Lemma}

\proof
Noting that $F^2$ is a polynomial in $s,\lambda,\phi$, we see that 
\begin{align*}
\frac{\partial F}{\partial \lambda} & = \frac{1}{2F^2} \frac{\partial F^2}{\partial \lambda} F=\tilde P_1 F,\\
\frac{\partial F}{\partial \phi} & = \frac{1}{2F^2} \frac{\partial F^2}{\partial \phi} F=\tilde P_2 F,\\
\frac{\partial F}{\partial s} & = \frac{1}{2F^2} \frac{\partial F^2}{\partial s} F=\tilde P_3 F,
\end{align*}
where $\tilde P_1,\tilde P_2, \tilde P_3$ are rational functions of $s,\lambda,\phi$. This means that the formal power series $F$ is holonomic of rank $1$. Arguing as in the previous lemma, we find $P_1,P_2 \in \R(s,\lambda)[\phi]$ of degree $\leq 3$  such that 
\begin{align*}
\frac{\partial F}{\partial \lambda}(s,\lambda,\phi(\lambda)) & = P_1(s,\lambda,\phi(\lambda)) F(s,\lambda,\phi(\lambda))\\
\frac{\partial F}{\partial \phi}(s,\lambda,\phi(\lambda)) & = P_2(s,\lambda,\phi(\lambda)) F(s,\lambda,\phi(\lambda)).
\end{align*}

By the chain rule we find that 
\begin{align*}
\frac{d}{d\lambda} F(s,\lambda,\phi(\lambda)) & = \frac{\partial F}{\partial \lambda}(s,\lambda,\phi(\lambda)) + \frac{\partial F}{\partial \phi}(s,\lambda,\phi(\lambda)) \phi'(\lambda)\\
& = \left[P_1(s,\lambda,\phi(\lambda))+P_0(\lambda,\phi(\lambda)) P_2(s,\lambda,\phi(\lambda))\right]F(s,\lambda,\phi(\lambda)).
\end{align*}
Letting $Q_1$ be the reduction of $P_1+P_0P_2$ modulo $P$ finishes the proof.
\endproof

Since $\phi(\lambda)$ is algebraic of rank $4$ and $F$ is holonomic of degree $1$, it follows that $F(s,\lambda,\phi(\lambda))$ is holonomic of degree $4$. More precisely, we obtain the following lemma.

\begin{Lemma}
There are $Q_i \in \R(s,\lambda)[\phi]$ of degree $\leq 3$ such that 
\begin{equation} \label{eq_derivative_f}
\frac{d^i F}{d \lambda^i}(s,\lambda,\phi(\lambda)) = Q_i(s,\lambda,\phi(\lambda)) F(s,\lambda,\phi(\lambda)).
\end{equation}
\end{Lemma}

\proof
We construct $Q_i$ by induction on $i$. With $Q_0:=1$ and with $Q_1$ as in the previous lemma, the cases $i=0,1$ are already done. 

Once $Q_i$ is defined, we take the derivative to find that
\begin{displaymath}
\frac{d^{i+1} F}{d \lambda^{i+1}}(s,\lambda,\phi(\lambda)) = \tilde Q_{i+1}(s,\lambda,\phi(\lambda)) F(s,\lambda,\phi(\lambda)),
\end{displaymath}
where 
\begin{align*}
\tilde Q_{i+1}(s,\lambda,\phi) & :=\frac{\partial Q_i}{\partial \lambda} (s,\lambda,\phi)+\frac{\partial Q_i}{\partial \phi} (s,\lambda,\phi) P_0(\lambda,\phi) \\
& \quad + Q_i(s,\lambda,\phi) Q_1(s,\lambda,\phi) \in \R(s,\lambda)[\phi].
\end{align*}
We let $Q_{i+1} \in \R(s,\lambda)[\phi]$ be the remainder of $\tilde Q_{i+1}$ by division  by $P$. Then \eqref{eq_derivative_f} is satisfied with $i$ replaced by $(i+1)$. 
\endproof

Write 
\begin{displaymath}
b(s,\lambda)=\sum_{l=0}^\infty \frac{b_l}{l!} \lambda^l, \quad b_l \in \R(s).
\end{displaymath}
The statement of the proposition is equivalent to $\deg b_l \leq 2l$ (in particular $b_l \in \R[s]$, which is a non trivial fact).

\begin{Lemma}
Suppose that $b_m$ is a polynomial in $s$ for all $m \leq l$. Then $(s-1)b_{l+1}$ is a polynomial in $s$.  
\end{Lemma}

\proof
The space of polynomials in $\phi$ of degree $\leq 3$ with coefficients in $\R(s,\lambda)$ is a vector space of dimension $4$ over the field $\R(s,\lambda)$. Hence there exist rational functions $\tilde R_i \in \R(s,\lambda)$ which are not all zero with 
\begin{displaymath}
\sum_{i=0}^4 \tilde R_i(s,\lambda) \frac{1}{i!} Q_i (s,\lambda,\phi)=0.
\end{displaymath}
Clearing denominators and multiplying by $F(s,\lambda,\phi(\lambda))$ we find polynomials $R_i \in \R[s,\lambda]$ with 
\begin{displaymath}
\sum_{i=0}^4 R_i(s,\lambda) \frac{1}{i!} \frac{d^i F}{d \lambda^i}(s,\lambda,\phi(\lambda)) =0.
\end{displaymath}

Write $R_i=\sum_{k=0}^\infty \frac{1}{k!} R_{ik}(s)\lambda^k$ with $R_{ik} \in \R[s]$. The explicit computation of these polynomials is tedious, but a computer algebra software can handle this very quickly. We only need the following properties. 
\begin{align} 
R_{ik} & = 0 \quad \text{ if }k<i-1,\label{eq_rvanishing_rs} \\
 R_{i,i-1} & = c_i(s-1), c_i>0, i=1,\ldots,4 \label{eq_special_values_r}.
\end{align}

We have 
\begin{displaymath}
0=R_0 (b(\lambda)-s)+\sum_{i=1}^4 R_i b^{(i)}(\lambda)=R_0 (b(\lambda)-s)+\sum_{i=1}^4 \sum_{j,k=0}^\infty \frac{R_{ik}}{k!} \frac{b_{j+i}}{j!} \lambda^{j+k}.
\end{displaymath}

Comparing the coefficient of $\lambda^l, l \geq 0$ and using \eqref{eq_rvanishing_rs} yields
\begin{displaymath}
\sum_{i=1}^4 \binom{l}{l+1-i} R_{i,i-1}b_{l+1}+ \sum_{m=0}^{l} \sum_{i=0}^4 \binom{l}{m-i} R_{i,l+i-m}b_m=R_{0,l}s.
\end{displaymath}

By \eqref{eq_special_values_r}, $\sum_{i=1}^4 \binom{l}{l+1-i} R_{i,i-1}$ is a non-zero scalar multiple of $(s-1)$. The second sum and the right hand side are polynomials in $s$ by assumption, hence $(s-1)b_{l+1}$ is a polynomial in $s$. 
\endproof

\begin{Lemma}
Suppose that $b_m$ is a polynomial in $s$ for all $m \leq l+1$. Then $(3s+1)b_{l+2}$ is a polynomial in $s$.
\end{Lemma}

\proof
We argue as in the previous lemma, but this time using that the polynomials $\frac{1}{i!} Q_i(s,\lambda,\phi),i=1,\ldots,5$ are linearly dependent. We then find polynomials $\hat R_i \in \R[s,\lambda],i=1,\ldots,5$ such that 
\begin{displaymath}
\sum_{i=1}^5 \hat R_i(s,\lambda) \frac{1}{i!} \frac{d^iF}{d\lambda^i}(s,\lambda,\phi(\lambda))=0.
\end{displaymath}

We decompose 
\begin{displaymath}
 \hat R_i=\sum_{k}^\infty \frac{1}{k!} \hat R_{ik}(s)\lambda^k.
\end{displaymath}

We only need that   
\begin{equation} \label{eq_bound_hat_ri}
\deg \hat R_i\leq 2(3-i), \quad  \hat R_{ik}=0 \text{ if } k<i-2, 
\end{equation}
and 
\begin{equation} \label{eq_explicit_hat_rij}
 \hat R_{i,i-2}=\hat c_i (3s+1), \hat c_i>0, i=2,\ldots,5.
\end{equation} 

Comparing coefficients as above gives us for every $l \geq 0$
\begin{equation} \label{eq_recursion_bs_partb}
\sum_{i=2}^5 \binom{l}{l+2-i} \hat R_{i,i-2} b_{l+2}+\sum_{m=0}^{l+1} \sum_{i=1}^5 \binom{l}{m-i} \hat R_{i,l+i-m} b_m=0.
\end{equation}

By \eqref{eq_explicit_hat_rij}, $\sum_{i=2}^5 \binom{l}{l+2-i} \hat R_{i,i-2}$ is a non-zero multiple of $(3s+1)$. The second sum is  a polynomial in $s$ by assumption and the statement follows. 
\endproof

\begin{Lemma}
 $b_{l}$ is a polynomial in $s$ of degree $\leq 2l$ for all $l$.
\end{Lemma}

\proof
We use induction on $l$. Clearly $b_0=1$ and $b_1=3s+1$ are polynomials. If all $b_m,m \leq l+1$ are polynomials and $l \geq 0$, then $(s-1)b_{l+2}$ and $(3s+1)b_{l+2}$ are polynomials. Since $(s-1)$ and $(3s+1)$ are relatively prime, this implies that $b_{l+2}$ is a polynomial in $s$. 

The induction hypothesis and \eqref{eq_bound_hat_ri} imply that for each $1 \leq i \leq 5, 0 \leq m \leq l+1$ we have 
\begin{displaymath}
\deg \hat R_{i,l+i-m}b_m = \deg \hat R_{i,l+i-m}+\deg b_m \leq 2(l-m+3)+2m = 2(l+3). 
\end{displaymath}
From \eqref{eq_recursion_bs_partb} and $\deg(3s+1)=2$ we conclude that $\deg b_{l+2} \leq 2(l+2)$.  
\endproof

\def\cprime{$'$}

\end{document}